\documentclass[12pt]{article}

\title{{ \bf Elimination Theory \\ in Codimension Two }}
\author{ Alicia Dickenstein and Bernd Sturmfels}
\date{}

\usepackage{amsthm}
\usepackage{amsmath}
\usepackage{amsfonts} 

\newcommand{\baseRing}[1]{\ensuremath{\mathbb{#1}}}
\newcommand{\Z}{\baseRing{Z}}
\newcommand{\R}{\baseRing{R}}
\newcommand{\Q}{\baseRing{Q}}
\newcommand{\C}{\baseRing{C}}

\newcommand{\Ch}{{\mathcal C}}
\newcommand{\res}{{\rm Res}}

\theoremstyle{plain}
\newtheorem{theorem}{Theorem}[section]
\newtheorem{lemma}[theorem]{Lemma}
\newtheorem{corollary}[theorem]{Corollary}
\newtheorem{proposition}[theorem]{Proposition}

\newtheorem{definition}[theorem]{Definition}
\newtheorem{remark}[theorem]{Remark}

\newtheorem{example}[theorem]{Example}

\numberwithin{equation}{section}

\begin{document}
\maketitle

\begin{abstract} \noindent
New formulas are given for Chow forms, discriminants and resultants
arising from (not necessarily normal) toric varieties of codimension~$2$.
The Newton polygon of the discriminant is determined exactly.
\end{abstract}

\vskip -.5cm

\section{Introduction}

\vskip -.2cm

Sparse elimination theory concerns the study
of Chow forms and discriminants associated with toric
varieties, that is, subvarieties of projective space
which are parametrized by monomials \cite{gkzbook}, \cite{StuSparse}.
This theory has its
origin in the work of Gel'fand, Kapranov and Zelevinsky
on multivariate hypergeometric functions \cite{gkz1}. The singularities
of these functions occur on the projectively dual hypersurfaces to
the torus orbit closures on the given toric variety $X$.
The singular locus of the hypergeometric system
is described by the full discriminant of $X$,
which is a natural specialization of the Chow form.

Classical hypergeometric functions in one variable arise
when $X$ is a toric hypersurface,
defined by one homogeneous binomial
equation $\, x_1^{b_1} \cdots x_r^{b_r} =
x_{r+1}^{b_{r+1}} \cdots x_{n}^{b_{n}} $.
The Chow form of this hypersurface $X$ is
just its defining polynomial.
The discriminant of $X$ equals, up to an integer factor,
\cite[\S 9.1]{gkzbook},
\begin{equation}
\label{introej1}
D_X \,\,\, = \,\,\, b_{r+1}^{b_{r+1}} \cdots b_{n}^{b_{n}} \cdot
x_1^{b_1} \cdots x_r^{b_r} - (-1)^{deg(X)}
b_1^{b_1} \cdots b_r^{b_r} \cdot
x_{r+1}^{b_{r+1}} \cdots x_{n}^{b_{n}} , \,\,\,
\end{equation}
and the full discriminant equals
$\, D_X \, $ times $\, \prod_{i=1}^n x_i^{deg(X)-b_i}$.
It is the purpose of this article to generalize these
formulas to toric varieties of
codimension $2$. 
Our motivations for this study include
hypergeometric functions \cite{rat}, 
Horn systems in two variables \cite{timur},
and their applications to theoretical physics
\cite{mp}.

We introduce our objects of study by means of an
example. Let $X$ be the toric $6$-fold
in projective $8$-space given parametrically
by the cubic monomials
$$ (a: b:\cdots :i) \, = \,
(
u_1 u_0^2 :
u_2 u_0^2 :
u_3 u_0^2 :
u_1 x^2 :
u_2 y^2 :
u_3 z^2 :
u_4 y z :
u_4 x z :
u_4 x y) . $$
The prime ideal of the toric variety $X$ is generated by
the $2 \times 2$-minors of
\begin{equation}
\label{introEx1}
\left( \begin{array}{ccc}
a & b & c \\
d g^2 & e h^2 & f i^2
\end{array} \right)
\end{equation}
Thus $X$ is arithmetically Cohen-Macaulay and has
degree $13$. The Chow form of $X$ is gotten
by eliminating the variable $t$ from
the $2 \times 2$-minors of
\begin{equation}
\label{introEx2}
\left(
\begin{array}{ccc}
a_0 + t a_1 & b_0 + t b_1 & c_0 + t c_1 \\
(d_0 + t d_1) (g_0 + t g_1)^2 &
(e_0 + t e_1) (h_0 + t h_1)^2 &
(f_0 + t f_1) (i_0 + t i_1)^2
\end{array} \right)
\end{equation}
The Chow form is an irreducible polynomial
of degree $26$ in the $18$ variables
$\,a_0,a_1,b_0,b_1,\ldots,i_0,i_1 \,$ having exactly
$57,726$ terms. It equals the
determinant
\begin{equation}
\label{introEx3}
\left(
\begin{array}{cccc}
123 &
124 &
125 &
126 \\
134 &
135 + 234 &
136 + 235 &
236 \\
135 &
136 + 145 + 235 &
146 + 236 + 245 &
246 \\
136 &
146 + 236 &
156 + 246 &
256
\end{array} \right)
\end{equation}
where $ijk$ is the $3 \times 3$-minor
with row indices $i$, $j$ and $k$ of the
$6 \times 3$-matrix
\begin{equation}
\label{introEx4}
\left(
\begin{array}{cccc}
a_0 & b_0 & c_0
\\
a_1 & b_1 & c_1
\\
d_0 g_0^2 & e_0 h_0^2 & f_0 i_0^2 \\
d_1g_0^2 + 2d_0 g_0 g_1 & e_1 h_0^2 + 2e_0 h_0 h_1 & f_1 i_0^2 + 2f_0
i_0 i_1\\
d_0g_1^2 + 2d_1 g_0 g_1 & e_0 h_1^2 + 2e_1 h_0 h_1 & f_0 i_1^2 + 2f_1
i_0 i_1\\
d_1 g_1^2 & e_1 h_1^2 & f_1 i_1^2
\end{array} \right)
\end{equation}
Note that the Chow form can also be written as a polynomial
of degree $13$ in the brackets
$\,[ab] = a_0 b_1 - a_1 b_0, \,
[ac] = a_0 c_1 - a_1 c_0, \, \ldots,\,
[hi] = h_0 i_1 - h_1 i_0$.
We obtain the full discriminant of $X$ from the Chow form by
substituting
\begin{equation}
\label{introEx5}
\left(
\begin{array}{cc}
a_0 & a_1 \\
b_0 & b_1 \\
\vdots & \vdots \\
i_0 & i_1
\end{array} \right) \quad \mapsto \quad
{\rm diag}(a,b,c,d,e,f,g,h,i) \cdot B,
\end{equation}
where $B$ is the $9 \times 2$-matrix with row vectors
$\,(1,0),(0,1),(-1,-1),(-1,0)$, $(0,-1),(1,1),(-2,0),(0,-2),(2,2)$.
The result of this substitution is the dual full discriminant $\tilde E_X$.
It has exactly twelve terms and factors as follows:
\begin{equation}
\label{Fourteen}
\tilde E_X \quad = \quad
2^{14} \cdot (a e h^2 - b d g^2 ) \cdot (a f i^2-c d g^2)
\cdot (b f i^2 - c e h^2 ) \cdot \tilde D_X,
\end{equation}
where the last factor $\tilde D_X$
is the irreducible polynomial
\begin{eqnarray*}
&
a^2 e^2 f^2 h^4 i^4
+ b^2 d^2 f^2 g^4 i^4
+ c^2 d^2 e^2 g^4 h^4 \\
&
-2 a b d e f^2 g^2 h^2 i^4
- 2 acdf e^2 g^2 h^4 i^2
-2 bcef d^2 g^4 h^2 i^2
\end{eqnarray*}
Replacing each variable in $\tilde D_X$
by its reciprocal, that is,
$\, a\mapsto 1/a, b \mapsto 1/b, \ldots \,$ and
clearing denominators, we get the discriminant $D_X$,
an irreducible polynomial of degree $10$ which defines the
hypersurface projectively dual to $X$.

In this paper we establish exact formulas for the
Chow form (Theorems \ref{chowth} and \ref{bezout}),
the full discriminant (Proposition \ref{fullpro}),
and the discriminant (Theorem \ref{algoth}) associated with an
arbitrary
toric variety $X $ of codimension $2$ in a projective space.
A combinatorial construction is given for the secondary polygon
(Theorem \ref{chpfdth}) and the
Newton polygon of the discriminant (Theorem \ref{nq}).
This construction shows that the
dual variety $X^\vee$ is a hypersurface if and only if
the secondary polygon is not centrally symmetric (Corollary \ref{D1}).
In Section 5 we study mixed resultants,
that is, we apply our theory to codimension $2$ toric varieties
which arise from the Cayley trick \cite[\S 3.2.D]{gkzbook}

The toric $6$-fold $X$ in our example does
arise from the Cayley trick. This can be seen from
the defining parametrization $\,( u_1 u_0^2 : \cdots : u_4 x y) $.
Hence the discriminant $D_X$ is actually a resultant.
Indeed, if we eliminate $x,y,z$ from
\begin{equation}
\label{ThreeBinoOneTrino}
a + d \cdot x^2 \,\, = \,\,
b + e \cdot y^2 \,\, = \,\,
c + f \cdot z^2 \,\, = \,\,
g \cdot y z + h \cdot x z + i \cdot xy \,\, = \,\, 0
\end{equation}
then the result is precisely the six-term discriminant $D_X$ described
above.

\medskip

\section{The Chow form} \label{chsection}

Let $B = (b_{i\ell})$ be an $n \times 2$-integer matrix of rank $2$
with
both column sums equal to zero. The {\it lattice
ideal} $I_B$ is the ideal in $k[x_1,\ldots,x_n]$,
$k$ a field of characteristic zero, generated by the
binomials $x^{u_+} - x^{u_-}$ where $u = u_+ - u_-$
runs over the two-dimensional lattice $L_B \subset \Z^n$ spanned by the
columns of $B$. The minimal
generators and the higher syzygies of $I_B$ are
described explicitly in \cite{ps}.
The ideal $I_B$ is homogeneous with respect to the usual $\Z$-grading
and hence defines a subscheme $X_B$ of projective space
${\bf P}^{n-1}$. The lattice
ideal $I_B$ is prime if and only if $\Z^n / L_B$ is a free abelian
group, or equivalently, if and only if the rows of
$B$ generate the
two-dimensional lattice $\Z^2$.

In this section we compute the Chow form and the Chow
polygon of the projective scheme $X_B$.
The {\it degree} of $X_B$, denoted
$\,d_B = degree(X_B)$, is the
number of intersection points with
a generic $2$-plane in ${\bf P}^{n-1}$.
Let $Y = (y_{i\ell})$ be an $n \times 2$-matrix of indeterminates.
It represents a generic parametric line
$(y_{11} + t y_{12}, \ldots, y_{n1} + t y_{n2})$
in ${\bf P}^{n-1}$. Following \cite[\S 3.2.B]{gkzbook},
the {\it Chow form} $\tilde{\Ch}_B$ of the homogeneous
lattice ideal $I_B$ is the unique (up to sign)
irreducible homogeneous polynomial in $\Z[y_{i\ell}]$ which
vanishes if and only if the corresponding
line in ${\bf P}^{n-1}$ meets $X_B$.
The degree of $\tilde{\Ch}_B$ equals
$\, 2 \cdot d_B $.

Classical invariant theory (cf.~\cite[Proposition 3.1.6]{gkzbook})
tells us that the Chow form $\tilde{\Ch}_B$ can be written
(non-uniquely)
as a polynomial of degree $d_B$ in the
{\it (dual) Pl\"ucker coordinates} of a
generic line, which we write as brackets
$$ [\, i \, j \, ] \quad := \quad
y_{i1} y_{j2} - y_{i2} y_{j1} \qquad
\hbox{for} \quad 1 \leq i < j \leq n . $$
We further introduce a non-negative integer $\nu_{ij}$
for any $\,1 \leq i < j \leq n \,$ as follows:
if the $i$-th row vector and the $j$-th row
vector of $B = (b_{i\ell})$ have the same sign in one of the two
coordinates
then set $\nu_{ij} = 0$; otherwise we set
\begin{equation}
\label{nu}
\nu_{ij} \quad := \quad
min \, \bigl\{
|b_{i1} b_{j2}|,
|b_{i2} b_{j1}| \bigr\}.
\end{equation}
Thus, $\nu_{ij} =0$ unless $b_i$ and $b_j$ lie in the interior
of opposite quadrants. Let
\begin{equation}
\label{polys}
H_\ell(t) \,\, = \,\, \prod_{i:b_{i\ell} > 0} (y_{i1} + y_{i2}
t)^{b_{i\ell}} -
\prod_{i:b_{i\ell} < 0} (y_{i1} + y_{i2} t )^{-b_{i\ell}}\,,
\,\,\,\,\ell
=1,2.
\end{equation}
We regard $H_1$ and $H_2$ as polynomials in a single variable $t$ with
coefficients
in $\Z[ y_{i\ell}\, , i=1,\ldots, n, \, \ell =1,2]$.
Let $\beta_\ell$ denote the sum of the positive
entries in the $\ell$-th column of $B$, for $\ell=1,2$. Clearly,
$degree( h_\ell) = \beta_\ell \, , \ell=1,2.$

\begin{theorem}\label{chowth}
The Chow form of the codimension $2$ lattice ideal $I_B$ equals
$$
\tilde{\Ch}_B \quad = \quad
\frac { \res_t \,( {\it H}_1, {\it H}_2) }
{\prod_{1 \leq r < s \leq n}
[\, r \, s \,]^{\nu_{rs}}},
$$
where $\res_t$ denotes the Sylvester resultant
of two univariate polynomials.
\end{theorem}

\begin{proof}
The binomials
$\prod_{b_{ij} > 0} x_i^{b_{ij}} - \prod_{b_{ij} < 0} x_i^{-b_{ij}}\, ,
\,
j=1,2,$
defined by the two columns of $B$
determine a complete intersection $Y_B$ of degree $\beta_1 \beta_2$
in ${\bf P}^{n-1}$ which coincides with $X_B$ over $\left(k^*\right)^{n-1}$.
The irreducible decomposition of $Y_B$
consists of the components of $X_B$ -- of which there is
only one if $\Z^n / L_B$ is free abelian --
together with subschemes supported on coordinate flats
$\,x_r = x_s = 0,$ whose Chow forms are
the bracket monomials $[\, r \, s \,].$
The theorem will be proved if we show that
the cycle $\,\{x_r = x_s = 0\}\, $ occurs with
multiplicity $\nu_{rs}$ in the complete intersection.

Suppose first that $\nu_{rs}=0$. We may assume that $b_{r1},b_{s1} \geq
0.$
Then, $\{x_r = x_s = 0\}$ is not contained in $Y_B$, and thus occurs
with
multiplicity $0$. Suppose now that $\nu_{rs} >0$.
We may assume that $b_{r1},b_{r2} >0$
and $b_{s1},b_{s2} < 0.$ Then, $\{x_r = x_s = 0\}$ is contained in
$Y_B$,
and
after localizing and changing variable names, we are
lead to the following situation:
let $a,b,c,d \in \Z_{>0}, \, ad> bc$ and $\alpha, \beta \not=0 $ in an
extension field $K$ of $k$, and consider the univariate resultant
$$r \,\,:= \,\,\res_t \bigl(
(x_0+x_1t)^a - \alpha (y_0 + y_1t)^b, (x_0+x_1t)^c - \beta (y_0 +
y_1t)^d
\bigr).$$
We want to show that $x_0 y_1 - y_0 x_1$ appears with exponent
$bc$ as a factor of $r$.

Indeed, when $x_1, y_1 \not=0,$ the condition $x_0 y_1 - y_0 x_1 =0$
holds if and only if there exists $t$ such that $x_0 + x_1t = y_0 + y_1
t
=0,$
and so $x_0 y_1 - y_0 x_1$ occurs in $r$ with exponent $\mu$ equal to
the intersection multiplicity at the origin of the artinian
ideal $\,I = \langle x^a - \alpha y^b, x^c - \beta y^d \rangle\,$
in $K[x,y]$. We claim $\mu = bc$.

The given equations are a Gr\"obner basis with leading terms
$x^a$ and $\beta y^d$, for the term order defined by
${\rm weight}(x) = b+d $ and ${\rm weight}(y) = a+c $.
Hence $\, {\rm dim}_K K[x,y] / I = a d$, that is,
there are $ad$ roots in the affine plane counting multiplicity.
Of those, $ad-bc$ lie in the torus, i.e., have both coordinates
non-zero. No root of $I$ has precisely one zero coordinate.
Therefore the multiplicity of $I$ at the origin is the difference
$\, \mu \, = \, ad \, - \, (ad - bc) \, = \, b c $.
\end{proof}

\begin{corollary} \label{degreeco}
The degree of a homogeneous
lattice ideal $I_B$ of codimension two
can be computed from the defining $n \times 2$-matrix $B$
by the following formula
$$ {\rm degree}(X_B) \quad = \quad
\beta_1 \beta_2 \, - \, \sum_{1 \leq r < s \leq n} \!\! \nu_{rs}. $$
\end{corollary}

The polynomial ring $\Z[y_{i\ell}]$ has a natural $\Z^n$-grading
defined by $\, deg(y_{i\ell}) = e_i$, the $i$-th unit vector.
The {\it Chow polytope} $CP_B$ is, by definition
\cite[\S 6.3]{gkzbook}, the
convex hull in $\R^n$ of the degrees of all monomials
appearing in the expansion of $\,\tilde{\Ch}_B$. Its faces correspond
to toric deformations of the algebraic cycle $X_B$.

We assume that the row vectors $b_1,b_2,\ldots,b_m$ of
the matrix $B$ are ordered counterclockwise in cyclic order, and that
$b_{m+1}, \ldots,b_n =0$. It may happen
that $b_{i+1}$ is a positive multiple of $b_{i}$.
Let $P_B$ denote the unique (up to translation) lattice
polygon whose boundary consists of the directed edges
$b_1,b_2,\ldots,b_m$. For each vector $b_i = (b_{i1}, b_{i2})$,
the linear functional
$$ u = (u_1,u_2) \quad \mapsto \quad \det(b_i,u) \,=
\, b_{i1} u_2 - b_{i2} u_1 $$
attains its {\sl minimum} value over $P_B$ at the edge parallel to
$b_i$
for $i=1,\ldots,m$ and is zero for $i=m+1,\ldots,n.$
Let $\mu_i$ denote the {\sl maximum} value of the linear functional
$u \mapsto \det(b_i,u)$ as $u$ ranges over the polygon $P_B$.
For $i=1,\ldots,m$, this maximum is attained at a unique vertex of
$P_B$ unless
$b_j = \lambda b_i$ for some $j$ and $\lambda < 0 $.
For every lattice point $v$ in $P_B$, the quantity
\begin{equation} \label{coords}
v^{(i)} \quad := \quad \mu_i - \, \det(b_i,v)
\end{equation}
is a non-negative integer, invariant under translation of $P_B$.
The vector $\, (
v^{(1)}, v^{(2)},
\ldots, v^{(n)})\,$ expresses the point
$v$ in $P_B$ in intrinsic coordinates.

\begin{theorem} \label{chpth}
The Chow polygon $CP_B$ of a codimension $2$
lattice ideal $I_B$ is the image of the polygon $P_B$ under
the affine isomorphism $v \mapsto (
v^{(1)} \! ,
\ldots,
v^{(n)})$.
\end{theorem}

The proof of this theorem will be given in the next section, after
Gale duality and duality of Pl\"ucker coordinates have been introduced.
See Theorem \ref{chpfdth} for the same theorem in dual formulation.
Theorems \ref{chpth} and \ref{chpfdth} will then derived
from the constructions in Sections 7.1.D and 8.3.B
of \cite{gkzbook}.

\begin{example} \rm
For the example in the Introduction we take
$\, b_1 = (1,0), \,
b_2 = (1,1), \,
b_3 = (2,2), \,
b_4 = (0,1), \,
b_5 = (-1,0), \,
b_6 = (-2,0), \,
b_7 = (-1,-1), \,
b_8 = (0, -1), \,
b_9 = (0,-2)\,$ and $P_B$ the hexagon with
vertices $\, (0,0), (1,0),
(4,3)$, $ (4,4),(1,4), (0,3)$. The edges of $P_B$
are labeled by the variables as follows:
$\, a,\{f,i\} ,b, \{d,g\} ,c, \{e,h\} $, and we have
$\mu = (4, 3, 6, 0, 0, 0, 1, 4, 8)$.
The twelve points on the boundary of $P_B$ correspond
to the twelve monomials in the expansion of $\tilde E_X$ .
For instance, the vertex $v = (0,0)$ has intrinsic coordinates
$\,(v^{(1)}, \ldots, v^{(9)}) = (4, 3, 6, 0, 0, 0, 1, 4, 8) $ and
corresponds to $a^4 c e^4 f^3 h^8 i^6$.
\end{example}

For any $v \in P_B$, the
coordinate sum $\,\sum_{i=1}^{n} v^{(i)} \,$ coincides with
$ \sum_{i=1}^{n} \mu_i$, and this equals
the degree of the Chow form
$\tilde{\Ch}_B$ as a polynomial in the $y_{i\ell}$. From this we get
an alternative formula for the degree of our lattice ideal.

\begin{corollary} \label{degmu}
The degree of the variety $X_B$ equals $ \,
d_B \, = \,
\frac{1}{2} \cdot \sum_{i=1}^n \mu_i $
\end{corollary}

Counting lattice points in the polygon $P_B$ gives an upper bound for
the number of monomials appearing in the full discriminant $D_X \,$
(see \S \ref{full} below):

\begin{remark}\label{Pick}
The number of  lattice points in the polygon $P_B$ equals
$$
1 \, + \,\frac{1}{2}
\biggl(
\sum_{i=1}^n gcd(b_{i1},b_{i2}) \,\, + \,\
\sum_{1 \leq i < j \leq n}
(b_{i2} b_{j1} -
b_{i1} b_{j2} )
\biggr)
$$
\end{remark}

\begin{proof}
This is a reformulation of {\sl Pick's formula}
which states that the area
of a lattice polygon equals the number of lattice points in that
polygon
minus half the number of lattice points in its boundary, minus one.
\end{proof}

If the lattice ideal $I_B$ is a complete intersection
then the denominator in Theorem \ref{chowth} is $1$
and we get a determinantal formula for the Chow form, namely,
$\,\tilde{\Ch}_B \,$ equals the univariate resultant
in the numerator, which can be computed
as the determinant of a Sylvester or B\'ezoutian matrix.

It would be desirable to have a division-free determinantal formulas
for the Chow form $\tilde{\Ch}_B$ of any codimension $2$ lattice ideal.
At this time we know such formulas only for special classes of matrices
$B$. We present a formula for a class which includes
the example in the Introduction. Recall
from \cite{ps} that the lattice ideal $I_B$ is Cohen-Macaulay
if and only if $I_B$ is generated by the $2 \times 2$-minors
of a $2 \times 3$-matrix of monomials in $\,x_1,\ldots,x_n$:
$$ \left( \begin{array}{ccc}
m_1 & m_2 & m_3 \\
m_4 & m_5 & m_6
\end{array} \right).$$
Let $d_i$ denote the total degree of the monomial $m_i$.
In order for the lattice ideal
$I_B$ to be homogeneous it is necessary and sufficient that
$$
d_1 + d_5 \, = \, d_2 + d_4
\quad \hbox{and} \quad
d_1 + d_6 \, = \, d_3 + d_4 . $$
For the following discussion we make an even
more restrictive assumption:
\begin{equation}
\label{RestAss}
d_1 = d_2 = d_3 \,\, \geq \,\, d_4 = d_5 = d_6.
\end{equation}
We introduce four new indeterminates $\,s,t, u,v$.
Let $\, m_i[t] \,$ denote the image of the monomial $\, m_i \,$
under the substitution
$\, x_i \,\mapsto \, y_{i1} + y_{i2} t \,$ for $i=1,2,\ldots,n$.
We define the {\it B\'ezout polynomial} to be the following expression:
$$ \frac{1}{(s-u)(t-v)} \cdot
{\rm det} \!\!
\left( \begin{array}{ccc}
\! m_1[t]+ m_4[t]\cdot s & \, m_1[t]+ m_4[t] \cdot u & \, m_1[v]+
m_4[v]
\cdot u \\
\! m_2[t]+ m_5[t]\cdot s & \, m_2[t]+ m_5[t] \cdot u & \, m_2[v]+
m_5[v]
\cdot u \\
\! m_3[t]+ m_6[t]\cdot s & \, m_3[t]+ m_6[t] \cdot u & \, m_3[v]+
m_6[v]
\cdot u \\
\end{array} \! \right)$$
Set $\delta := d_1 + d_4 $. The B\'ezout polynomial can be
written uniquely in the form
$$
(1,v,v^2\ldots, v^{d_1-1},\, u,u v, uv^2,\ldots, uv^{d_4-1}) \cdot
{\bf B} \cdot
\left( \begin{array}{c} 1 \\ t \\ \vdots \\ t^{\delta-1}
\end{array} \right) , $$
where $ {\bf B} = {\bf B}(y_{ij})$ is a certain $\delta \times \delta$-matrix
with entries in $k[y_{11},y_{12},\ldots , y_{n2}]$.

\begin{theorem} \label{bezout}
If $I_B$ is a Cohen-Macaulay lattice ideal of codimension $2$
satisfying (\ref{RestAss}) then its
Chow form $\tilde{\Ch}_B$ equals the determinant of
$\,{\bf B}(y_{ij})$.
\end{theorem}

\begin{proof} Consider the rational normal scroll of type $(d_1,d_4)$,
a toric surface of degree $\delta$ in a projective space
of dimension $\delta+1$. Its Chow form
has an exact determinantal
formula in terms of a B\'ezout matrix.
A nice proof of this fact follows from recent results of
Eisenbud and Schreyer \cite{es}, since the rational normal scroll
is given by the $2 \times 2$-minors of a matrix of variables.
This Chow form is the unmixed, sparse resultant for three
polynomials with support
$$ \{
1,t,t^2, \ldots,t^{d_1}, \,
s,st,st^2, \ldots,st^{d_4} \}. $$
The three polynomials
$\, m_i[t] + m_{i+3}[t] \cdot s \,$ have exactly this support.
Our formula is gotten by specializing the
B\'ezout matrix for the scroll.
\end{proof}

\begin{example} \rm
The ideal in (\ref{introEx1}) satisfies the
hypotheses of Theorem \ref{bezout}, with $\delta = 4$.
The matrix (\ref{introEx3}) is precisely the
matrix $\,{\bf B}(y_{ij})$ in this case. \qed
\end{example}

\section {The full discriminant} \label{full}

There are two different ways of presenting a toric variety
of codimension two: by an $n \times 2$-matrix $B$ as in \cite{ps}, or
by an $(n-2) \times n$-matrix $A$ as in \cite[\S 5.1]{gkzbook}.
The two matrices are {\it Gale dual}, which means
that the image of $B$ equals the kernel of $A$.
Up to this point in the paper, we have only used
the $B$-representation. We now make a switch and introduce
the $A$-representation.

Let $ A = (a_1,\ldots,a_n)$ be an $(n-2)\times n$-integer matrix
of rank $n-2$, and suppose there exists a vector
$w \in \Q^{n-2}$ such that $w \cdot a_i = 1$
for $i=1,2,\ldots,n$. We can choose an integral
$n \times 2$ matrix $B$ whose columns are a $\Z$-basis of
$\ker_\Z(A)$. The matrix $B$ has rank $2$ and $A\cdot B = 0$.
It is unique modulo right multiplication by $GL(2,\Z)$.
Let $I_A = I_B$ denote the corresponding toric ideal in
$k[x_1,\ldots,x_n]$ and $X = X_A = X_B$ the corresponding
toric variety in ${\bf P}^{n-1}$.

Here it is important to note that not all integer matrices $B$ arise
as the Gale dual of some matrix $A$ as above. For this it is necessary
and
sufficient that $\Z^n/im_\Z(B)$ is torsion-free, or equivalently,
that the ideal $I_B$ is prime.

The {\it $A$-discriminant} $D_A$ is an irreducible polynomial
in $\Z[x_1,\ldots\!, x_n]$ which vanishes under a specialization
if the corresponding Laurent polynomial
$$ f \quad = \quad \sum_{i=1}^n \, x_i \cdot
t_1^{a_{i1}} t_2^{a_{i2}} \cdots t_{n-2}^{a_{i,{n-2}}}
\quad \hbox{where } x_1,\ldots,x_n \in \C^* $$
has a multiple root $(t_1,\ldots,t_{n-2})$ in $(\C^*)^{n-2}$.
Equivalently, the hypersurface $\{D_A =0\}$ is projectively dual
to the toric variety $X$, when the dual variety $X^\vee$ is a
hypersurface, and $D_A = 1$ otherwise; see \cite[\S 1.1 and \S
9.1]{gkzbook}.

In the next section we give a formula for the $A$-discriminant $D_A$
and its degree. In this section, we study a larger polynomial $E_A$
which contains $D_A$ as a factor. It is called the
{\it principal $A$-determinant} in \cite{gkzbook} but we prefer
the term {\it full discriminant}. Actually, our full discriminant
agrees with expression (1.1) in \cite[10.1.A]{gkzbook}, but there
is a slight inaccuracy in \cite[Theorem 10.1.2]{gkzbook} since
$E_A$ does not generally have content $1$. An extra integer factor
is needed. This integer factor would be $2^{14}$ for the example
(\ref{Fourteen}) in the Introduction.

Before stating the definition of $E_A$, we first review
the duality between primal and dual Pl\"ucker coordinates,
and see how it ties in with Gale duality.
For $ 1 \leq i < j \leq n$, let $B(i,j)$
the submatrix of $B$ consisting of the $i$-th and $j$-th rows,
and let $A \langle i,{j} \rangle$ denote the submatrix of $A$ gotten by
omitting the $i$-th and $j$-th columns. Here signs are adjusted
so that ${\det A \langle{i},{j} \rangle} = { \det B(i,j)}$, up to a global
constant.
In Section 2 we used an
$n \times 2$ matrix $Y= (y_{i\ell})$ of indeterminates.
The {\it dual Pl\"ucker coordinates}
of a line in ${\bf P}^{n-1}$ are
\begin{equation}
\label{bbrracc}
[\, i \, j \, ] \quad := \quad
{\rm det} \, Y(i,j) \,\,\, = \,\,\,
y_{i1} y_{j2} - y_{i2} y_{j1} \quad
\hbox{for} \quad 1 \leq i < j \leq n .
\end{equation}
Here we consider an $(n-2) \times n$-matrix
$Z=(z_{ij})$ of indeterminates.
The {\it primal Pl\"ucker coordinates}
of our line are the
$(n-2) \times (n-2)$-subdeterminants
$$ \langle\, {i}\, {j} \,\rangle \,\, = \,\, \det Z \langle {i}, {j}
\rangle
\quad \qquad \hbox{(with the sign adjusted as usual).} $$
The dual Chow form $\tilde{\Ch}_B$ is a polynomial of degree $d_B$
in the brackets (\ref{bbrracc}). Replacing
$\, [\, i \, j \, ] \, \mapsto \langle \, i \, j \, \rangle \,$
in $\tilde{\Ch}_B$ gives a homogeneous polynomial
of degree $(n-2) d_B$ in the variables $z_{ij}$. It is
denoted $\Ch_A$ and called the {\it primal Chow form}. Note that
$\Ch_A$ coincides with the {\it $A$-resultant}
defined in \cite[\S 8.2.A]{gkzbook}.

\begin{definition} \label{defFullDis}
\rm
The {\it full discriminant} $E_A$ is the image of the
primal Chow form $\Ch_A$ under the specialization
$ z_{ij} \mapsto a_{ij} x_j \,$ for
$i=1,\ldots,n-2$, $j=1,\ldots,n$.
\end{definition}

We next show how to compute the full discriminant
directly from the dual Chow form $\tilde{\Ch}_B$ and hence
from the formulas in Theorems \ref{chowth} and \ref{bezout}.

\begin{proposition}
\label{fullpro}
The full discriminant $E_A$ and the
dual Chow form $\tilde{\Ch}_B(y_{i\ell})$ are related
by the following formula:
\begin{equation}
\label{propFullPro}
E_A(x_1,\ldots,x_n) \,\,\, =\,\,\, (x_1 \dots x_n)^{d_B} \cdot
\tilde{\Ch}_B( b_{i\ell}/ x_i, \,i=1,\ldots,n , \ell=1,2).
\end{equation}
\end{proposition}

The exponent $d_B$ is the degree of the toric
variety $X$ and hence coincides with the normalized
volume of the $(n -3)$-dimensional polytope
$\,{\rm conv} (A)$.
Gale dual formulas for this volume are given in
Corollaries \ref{degreeco} and \ref{degmu}.

\begin{proof}
The specialization $\, z_{ij} \mapsto a_{ij} x_j \,$
in Definition \ref{defFullDis} is equivalent to
\begin{equation}
\label{AijSpec}
\langle \,{r}\, {s} \,\rangle \,\,\, \rightarrow \,\,\, \det
A \langle {r},{s} \rangle
\prod_{k \not= r,s} x_k \qquad \hbox{for } \, 1 \leq r < s \leq n
\end{equation}
at the level of primal Pl\"ucker coordinates.
The dual Chow form $\tilde{\Ch}_B$ is a $\Z$-linear
combination of bracket terms $\, \prod [\, r \, \, s \,]\,$
of degree $d_B$. If we substitute
$\, b_{i\ell}/ x_i \,$ for $\,y_{i\ell} \,$ in the expansion
of such a bracket term $\, \prod [\, r \, \, s \,]\,$ then we get
\begin{eqnarray*}
\prod [\, r \, \, s \,] \!\!\!\!\!\!\!\!\!\! & \longrightarrow \quad
\prod \, \bigl( {\rm det} B(r,s) / (x_r x_s) \bigr) \, = \,\,
\prod \, \bigl( {\rm det} A \langle r,s \rangle / (x_r x_s) \bigr) \,\, = \\
& \! (x_1 \cdots x_n)^{-d_B} \cdot
\prod \bigl( {\rm det} A \langle r,s \rangle \prod_{k \not= r,s} x_k \bigr)
\,\, \longleftarrow \,\, (x_1 \cdots x_n)^{-d_B} \cdot
\prod \langle\, r \, \, s \,\rangle
\end{eqnarray*}
Hence the specialized dual Chow form on
right hand side of (\ref{propFullPro}) equals the
specialization of the primal Chow form $\Ch_A$ under (\ref{AijSpec}),
as desired.
\end{proof}

It is known from \cite[Theorem 10.1.2]{gkzbook}
that the full discriminant $E_A$ is a product of
irreducible factors $D_{A'}$ where $A'$ ranges
over facial discriminants. In particular, each monomial
$x_i$ corresponding to a vertex $a_i$ of ${\rm conv}(A)$
appears to some positive power in the factorization of $E_A$.
It is curious to note that the monomial factors
disappear when we pass to dual coordinates. We define
the {\it dual full discriminant} by
specializing the dual Chow form:
\begin{equation} \label{fakech}
\tilde{E}_B(x_1,\dots,x_n)
\quad = \quad \tilde{\Ch}_B( \,b_{i\ell}\cdot x_i \, ,
\,\, i=1,\ldots,n \ , \ell=1,2).
\end{equation}
Proposition \ref{fullpro} is equivalent to the reciprocity formula:
\begin{equation}\label{fakeeq}
\tilde{E}_B(x_1,\ldots,x_n) \quad =
\quad (x_1 \dots x_n)^{d_B} \cdot E_A (1/x_1,\ldots, 1/x_n).
\end{equation}

\begin{lemma} \label{nomono}
The dual full discriminant $\tilde{E}_B$ has no monomial factors.
\end{lemma}

\begin{proof}
Suppose that the variable $x_i$ divides $\tilde{E}_B$.
Then every bracket monomial appearing in the dual Chow form
$\tilde{\Ch}_B$ contains the letter $i$. Equivalently, every
bracket monomial in the primal Chow form $\Ch_A$
contains a bracket $\, \langle \, r \, s \, \rangle \,$
with $r = i$ or $s=i$. In view of \cite[Theorem 8.3.3]{gkzbook},
this means that every regular triangulation of $A$
contains a simplex for which $a_i$ is not a vertex.
But this is false, since $a_i$ lies in every
maximal simplex of the {\it reverse lexicographic
triangulation} of $A$, for $x_i$ smallest; see
\cite[Proposition 8.6]{Stubook}.
\end{proof}

The {\it secondary polygon} $\Sigma(A)$ of the configuration $A$
coincides with the Newton polygon of the full discriminant $E_A$,
by \cite[Theorem 10.1.4]{gkzbook}. It is a
$2$-dimensional convex polytope lying in $\R^n$.
Let $P_B$ be the polygon
considered in Section \ref{chsection}. For $v \in P_B$, let
$(v^{(1)}, \ldots, v^{(n)})$ be the vector defined in (\ref{coords}).

\begin{theorem} \label{chpfdth}
The secondary polytope $\Sigma(A)$
is the image of the polygon $P_B$ under the affine isomorphism
which sends $\, v \, $ to
$ \,(d_B - v^{(1)}, \ldots, d_B - v^{(n)}).$
\end{theorem}

\begin{proof}
It suffices to prove this theorem for the case when all $b_i$
are non-zero. Indeed, if $b_{m+1} = \cdots = b_{n} = 0$ then
\cite[Theorem 10.1.2]{gkzbook} implies that
$$ E_A (x_1,\ldots,x_n) \quad = \quad \left(x_{m+1} \ldots
x_n\right)^{d_B}
\cdot E_{A'}(x_1,\ldots,x_m),$$
where $A'$ is a Gale dual of the
configuration $\,(b_1,\ldots, b_m)$.
Our assertion for $\Sigma(A')$ implies that for $\Sigma(A)$.
We hence assume that $b_i \not= 0$ for all $i =1,\ldots, n$.

Each vertex $w = (w_1,\ldots,w_n)$ of $\Sigma(A)$ corresponds uniquely
to the Gale dual of a regular triangulation $\Delta_w$, and hence to a
pair
of adjacent linearly independent vectors $b_k,b_{k+1}$
(indices are understood modulo $n$; recall that $\sum_i b_i =0$).
By \cite[Definition 7.1.6]{gkzbook},
the $i$-th coordinate of $w$ equals the sum of the
normalized volumes of those simplices in $\Delta_w$ which contain
the point $a_i$.
By Gale duality, $ w_i = \sum_{r,s} | \det(b_r,b_s)|$
where the sum is over all indices $r \not=i ,s \not=i$ such that
$b_k$ and $b_{k+1}$ lie in the cone spanned by $b_r$ and $b_s$.
Let $v_w$ be the vertex of $P_B$ between the edges
parallel to $b_k$ and $b_{k+1}$. We claim that
$v_w \in \Z^2$ is mapped to $w \in \Z^n$ under the
affine isomorphism given above.

We note that the maximum $\mu_i$ of the values $\det(b_i,v)$
is attained at the vertex $v \in P_B$ between
the edges parallel to two independent vectors $b_\ell, b_{\ell + 1}$
such
that
$\,\det(b_i, b_\ell) \geq 0$ and $\det(b_i, b_{\ell + 1}) <0 \,$
(indices modulo $n$). What we are claiming is the identity
$\, w_i \, = \, d_B - \det(b_i,v) + \det(b_i,v_w)$.

Let $C_k$ denote the set of index pairs $(r,s)$ such that
$b_k$ and $b_{k+1}$ lie in the cone spanned by $b_r$ and $b_s$.
The set $C_k$ is Gale dual to our regular triangulation, and,
hence $d_B = {\rm vol}({\rm conv}(A))$ equals
$\,\sum_{(r,s) \in C_k} | \det(b_r,b_s)|$. If we start drawing $P_B$
from the origin, then,
$v = \sum_{j=1}^{\ell} b_j$ and $v_w = \sum_{j=1}^{k} b_j.$
Our assertion takes the following form:
$$
\sum_{(r,s) \in C_k ,\,r \not=i ,s \not=i}
\!\!\!\!\!\!\!\!
| \det(b_r,b_s)| \,\,\,\, =
\sum_{(r,s) \in C_k} \!\! | \det(b_r,b_s)| - \sum_{j=1}^{\ell}
\det(b_i,b_j)
+
\sum_{j=1}^{k} \det(b_i,b_j).
$$
After erasing equal terms on both sides, the following remains to be
proved:
$$
\sum_{j=1}^{\ell}\det(b_i,b_j) - \sum_{j=1}^{k} \det(b_i,b_j) \,\, =
\sum_{j: (i,j) \in C_k} \!\!\! | \det(b_i,b_j)|.
$$
The proof is straightforward by a case distinction involving
the relative positions of
the vectors $b_i, b_k, b_{k+1} $ and $b_{\ell}$ in the plane.
\end{proof}

\vskip .1cm

\noindent {\sl Proof of Theorem \ref{chpth}: }
If $b_1,\ldots,b_n$ span the lattice $\Z^2$
then we find a corresponding matrix $A$, and Theorem \ref{chpth}
follows directly from Theorem \ref{chpfdth}
and the reciprocity formula (\ref{fakeeq}).
Otherwise, the scheme $X_B$ is the equidimensional union
of $r > 1$ torus translates of a fixed toric variety $X_{B'}$.
Following \cite[\S 4.1.A]{gkzbook}, the Chow form $\tilde{\Ch}_B$
factors into $r$ irreducible polynomials, each of which is
a torus translate of the irreducible Chow form $\tilde{\Ch}_{B'}$
Therefore the Chow polygon $CP_B$ equals $r \cdot CP_{B'}$.
The configuration $B'$ is $GL(\R^2)$-equivalent to $B$, and it
does possess a Gale dual $A'$. Our assertion holds for
$CP_{B'}$ and it follows for $CP_B$ by scaling.
\qed

\vskip .2cm

Let us now take a look at what happens to the formula
in Theorem \ref{chowth} under the specialization
$\,y_{i\ell} \mapsto b_{i\ell} \cdot x_i \,$ in (\ref{fakech}).
A line through the origin in $\R^2$ is said to be {\it relevant\/}
if it contains two vectors $b_r,b_s$ in opposite directions.
So, if the rows of $B$ are in general position, then there
are no relevant lines. The example in the introduction
has three relevant lines.

Consider the specializations of the two polynomials $H_\ell(t)$ in
(\ref{polys}):
\begin{equation}
\label{SpecialPolys}
h_\ell(t) \,\,\, := \,\,\, \prod_{i:b_{i\ell} > 0} (b_{i1} + b_{i2}
t)^{b_{i\ell}}
x_i^{b_{i\ell}} \, -\,
\prod_{i:b_{i\ell} < 0} (b_{i1} + b_{i2} t )^{-b_{i\ell}}
x_i^{-b_{i\ell}}
\,, \,\,\,\,\ell =1,2.
\end{equation}

\begin{remark}
\label{common}
The polynomials $h_1,h_2$ have a common factor
if and only if there exist a relevant line which is not a
coordinate axis.
\end{remark}

The presence of two vectors $b_r,b_s$ in opposite directions in the
interior of two quadrants then
causes the resultant $\res_t \,( h_1, h_2)$ to vanish. Also,
$\det(B(r,s))=0,$ while $\nu_{rs} \not=0.$ When there are
two opposite vectors on a coordinate axis, both numbers
are zero and $\det(B(r,s))^{\nu_{rs}}=1$. We deduce:

\begin{proposition} \label{fakepro}
Assume there are no relevant lines for the configuration $B$
except for the coordinate axes. Then the dual full discriminant equals
$$\tilde{E}_B \quad = \quad
\frac { \res_t \,( h_1, h_2) }
{\prod_{1 \leq r < s \leq n}
\det(B(r,s))^{\nu_{rs}}
{\prod_{1 \leq r < s \leq n}
(x_r \cdot x_s)^{\nu_{rs}}} }.
$$
\end{proposition}

In the next section we will show how to
use Theorem \ref{chowth} to compute discriminants
even if the hypothesis of the above proposition is not satisfied.

\section{The $A$-discriminant}

Let $A \in \Z^{(n-2)\times n}$ and $B \in \Z^{n \times 2}$ be
Gale dual matrices as before, and let $X $ be the corresponding
toric variety of codimension $2$ in $ {\bf P}^{n-1}$. The
{\it $A$-discriminant} $D_A$ is the defining
irreducible polynomial of the
dual variety $X^\vee$, unless $codim(X^\vee) > 1$ in which case
$\,D_A = 1$. Gel'fand, Kapranov and Zelevinsky
\cite[Theorem 10.1.2]{gkzbook} proved that
$D_A$ appears with exponent one in the factorization
of the full discriminant $E_A$. In this section we compute
$D_A$ and all other factors of $E_A$ in terms of the row
vectors $b_i \in \R^2$ of $B$.

Throughout this section we shall assume that
$b_i \not= 0$ for $i=1,\ldots,n$. This means
that $X$ is not a cone over a coordinate point,
or that $X^\vee$ does not lie in a
coordinate hyperplane.
All results in Section 4 require this hypothesis.

Each relevant line in the plane is identified with one of the
two primitive vectors $v \in \Z^2$ on that line. We abbreviate
$\, b^{(v)}_i \, := \, {\rm det}(b_i,v) $. With each such line $v$,
we associate a codimension one discriminant as in (\ref{introej1}).

\begin{equation}
\label{DDDv}
D_v \quad := \,\,\,
\prod_{j : \, b^{(v)}_j < 0 } \!\!\! (b^{(v)}_j)^{- b^{(v)}_j}
\prod_{i : \, b^{(v)}_i > 0 } \!\! x_i^{b^{(v)}_i} \,\, - \,\,
\prod_{i : \, b^{(v)}_i > 0 }\!\! (b^{(v)}_i)^{b^{(v)}_i}
\prod_{j : \, b^{(v)}_j < 0 } \!\!\! x_j^{- b^{(v)}_j}
\,\,\,
\end{equation}
Let $b_{i_1},\ldots,b_{i_s}$ be all the row vectors of $B$
which lie on the relevant line $v$. There is a unique
integer vector $(\lambda_1,\ldots,\lambda_s)$ such that
$\, b_{i_j} \, = \, \lambda_j \cdot v \,$ for $j=1,\ldots,s$.
We direct the primitive vector $v \in \Z^2$ so that the coordinate sum
$\,\alpha_v \, := \, \lambda_1 + \cdots + \lambda_s \,$
is nonnegative, and we define
$\, \delta_v \, := \,
\sum \{- \lambda_i : \lambda_i < 0 \} $.
Using this notation, Remark \ref{common} can now be refined as
follows:

\begin{remark}
\label{whichpower}
If $v = (v_1,v_2)$ is a relevant line for $B$ then
$v_1 + v_2 t$ appears with exponent
$\,\delta_v \cdot v_i \,$ in the factorization
of the polynomial $h_i(t)$ in (\ref{SpecialPolys}).
\end{remark}

Denote by $p_1(t),p_2(t)$ the respective remaining factors, that is,
\begin{equation} \label{pes}
h_\ell(t) \quad = \quad p_\ell(t) \,\cdot \!\!\!
\prod_{v \,{\text {relevant}}} \!\!\! (v_{1} + v_{2} t)^{\delta_{v}
\cdot v_\ell}
\quad \quad \ell=1,2.
\end{equation}
Now the resultant $\,r_B := \res_t (p_1,p_2) \,$
is a non-zero polynomial in $x_1,\ldots ,x_n$.
It is customary to call $r_B$ the
{\it residual resultant} of $h_1$ and $h_2$.
We shall prove the following formulas
for the full discriminant and the $A$-discriminant.

\begin{theorem}\label{algoth}
There exist monomials $x^u,x^{u'}$ and integers $\nu,\nu'$ such that
\begin{eqnarray*}
D_A (x_1,\ldots,x_n) & = \quad
(1/\nu) \cdot
x^u \cdot r_B(1/x_1,\ldots,1/x_n) \qquad \qquad \hbox{and} \\
E_A(x_1,\dots,x_n) & = \quad \nu' \cdot x^{u'} \cdot
D_A(x_1,\ldots,x_n) \cdot
\prod_{v \,{\rm {relevant}}} D_v(x_1,\ldots,x_n)^{\delta_v}
\end{eqnarray*}
\end{theorem}

\begin{proof}
We shall first prove the following claim about the full discriminant:
$$
r_B(1/x_1,\ldots,1/x_n) \cdot \!\!\!
\prod_{v \,{\rm {relevant}}} D_v(x_1,\ldots,x_n)^{\delta_v}
\quad \hbox{divides} \quad
E_A(x_1,\dots,x_n)
$$
in the Laurent polynomial ring $k[x_1,\ldots,x_n,
x_1^{-1},\ldots,x_n^{-1}]$.

Fix any relevant line $v$. Choose an isomorphism
in $SL_2(\Z)$ which maps $v$ to $(0,1)$, and apply
this isomorphism to the rows of $B$. Also reorder the
rows of $B$ so that the multiples of $v$ come first.
After this transformation, the first column of $B$ has the entries
$\,0, \ldots, 0, b_{s+1}^{(v)},b_{s+2}^{(v)}, \ldots,b_n^{(v)}$.

For $\ell=1,2$ and $i=1,\ldots,s$ {\bf only},
substitute $y_{i\ell} = b_{i\ell}/x_i$ into the Chow form $\Ch_B$.
Let $\tilde{H}_\ell$ be the polynomials resulting
from $H_\ell$ in (\ref{polys}) under the same substitution.
Then $\tilde{H}_1 = H_1$, but $\tilde{H}_2$
is divisible by $t^{\delta_v}$, and this is the
highest possible power of $t$ with this property
(cf.~Remark \ref{whichpower}).
Theorem \ref{chowth} implies that the specialized
Chow form factors, and one of its factors is
\begin{equation}
\label{aNiceIdentity}
\res_t(H_1, t^{\delta_v}) \quad = \quad \left( H_1(0)
\,\right)^{\delta_v}
\end{equation}
For all subsequent specializations, the Chow form factors
accordingly. When we substitute
$y_{i\ell} = b_{i\ell}/x_i$ for $i=s+1,\ldots,n, \ \ell=1,2,$
into $H_1(0)$ then we get the binomial $\,D_v \,$
in (\ref{DDDv}). Clearly, the residual resultant $r_B$ divides the
full discriminant $\tilde{E}_B$. The above claim follows from this.
Moreover, our argument shows that $D_v^{\delta_v}$ is the highest
power of $D_v$ which divides $E_A$.

Consider now the factorization formula given by
Gel'fand, Kapranov and Zelevinsky
in \cite[Theorem 10.1.2]{gkzbook}.
Under Gale duality,
the proper faces of the
polytope ${\rm conv}(A)$ which are not simplices correspond
to relevant lines $v$, and their
face discriminants are precisely the binomials $D_v$.
In other words, the full discriminant $E_A$ equals
the $A$-discriminant $D_A$ times the product of the
expressions $\,D_v^{\delta_v}\,$ where $v$ ranges over all
relevant lines. We conclude from our claim that
$\, r_B(1/x_1,\ldots,1/x_n) \,$ divides
$\,D_A(x_1,\ldots,x_n)\,$ in the Laurent polynomial ring.
Since $D_A$ is irreducible, both of our assertions follow.
\end{proof}

\vskip .1cm

We next compute the Newton polygon of the $A$-discriminant.
Define
$$ b_v \quad := \quad \alpha_v \cdot v \quad
= \quad b_{i_1} + \cdots + b_{i_s} $$
for any relevant line $v$.
It may happen that $b_v =0.$ We take all-non zero vectors $b_v$ and
all
vectors $b_i$ which do not lie in relevant lines,
and we order them counterclockwise in cyclic order.
Let $Q_B$
denote the unique (up to translation) lattice polygon whose
boundary consists of these directed edges. For any $i=1,\ldots,n$
and any lattice point $v$ in $Q_B$, we define
\begin{equation} \label{min}
\nu_i \, := \, \min \{ \, \det(b_i,u) \, , u \in Q_B \, \}
\quad \hbox{and} \quad
\overline{v}^{(i)} \, := \, \det(b_i,v) - \nu_i.
\end{equation}
Hence, $ \overline{v}^{(i)} \in \Z_{\geq 0}$
is the normalized lattice distance from $v$ to the boundary of $Q_B$,
in
the
direction orthogonal to $b_i.$

\begin{theorem} \label{nq}
The Newton polygon $N(D_A)$ of the $A$-discriminant $D_A$
is the image of the polygon
$Q_B$ under the affine isomorphism $v \mapsto \, ( \overline{v}^{(1)},
\ldots, \overline{v}^{(n)}).$
\end{theorem}

\begin{proof}
Suppose first that there are no relevant lines. Then, $Q_B = P_B$,
and the secondary polygon $\Sigma(A)$ and the Newton polygon $N(D_A)$
are equal up to translation.
More precisely, $\Sigma(A) = N(D_A) + \alpha$ where
$\alpha_i$ is the exponent of $x_i$ as a factor of $E_A$. Using
\cite[Theorem 10.1.2]{gkzbook} and Gale duality, we find
$$ \alpha_i \quad = \quad d_B \,\,\,\,\, - \,
\sum_{j : \,\det(b_i,b_j) > 0} \!\! \det(b_i,b_j). $$
In light of Theorem \ref{chpfdth}, it suffices to show that
$$ d_B - \mu_i + \det(b_i,v) \quad = \quad
\alpha_i + \det(b_i,v) - \nu_i
\qquad \hbox{for all} \,\, v \in Q_B \,, \, i = 1,\ldots,n. $$
After cancelling terms common to both sides, what remains to be shown
is
$$ \sum_{j : \,\det(b_i,b_j) > 0} \!\!\!\!
\det(b_i,b_j) \quad = \quad \mu_i - \nu_i. $$
This identity holds because both sides are equal to the
normalized lattice width of the polygon $Q_B = P_B$
in the direction orthogonal to $b_i$.

We next assume that relevant lines exist. Then $\nu_i$ generally
differs from $\nu'_i := \min \{ \det(b_i,u) \, , u \in P_B \} $.
The secondary polytope $\Sigma(A)$
equals $N(D_A) + \alpha$ plus the Minkowski sum of the Newton segments
of the binomials (\ref{DDDv}) where $v$ runs over all relevant lines.
Hence, up to lattice translation,
\begin{equation}
\label{polyPlusSeg}
P_B \quad = \quad Q_B \,\,\, + \, \sum_{v \,\, {\rm relevant}}
{\rm conv}\{\, 0, \, v \, \}.
\end{equation}
The minimum value of the linear functional $\det(b_i,*)$ over the
line segment $\,{\rm conv}\{\, 0, \, v \, \}\,$
is $ \det(b_i,v),$ when this value is negative and zero otherwise.
Therefore (\ref{polyPlusSeg}) translates into the identity
$$\nu_i' \quad = \quad \nu_i \,\,\,\, +
\sum_{v \, {\rm relevant}}\!\!
\delta_v \cdot {\rm min}\{ 0, \det(b_i,v)\} \qquad
\hbox{for} \,\, i=1,\ldots,n.$$
The argument for the case of
no relevant lines now completes the proof.
\end{proof}

We deduce the following formula for the degree of the $A$-discriminant:

\begin{corollary}\label{degdiscco}
$$ {\rm degree}(D_A) \quad = \quad
- \, \sum_{i=1}^n \nu_i $$
\end{corollary}

We can also extract the following characterization from Theorem
\ref{nq}.

\begin{corollary} \label{D1} The $A$-discriminant
$D_A$ is equal to $1$ if and only if the polygon $P_B$ is centrally
symmetric.
\end{corollary}

\begin{proof}
The condition $D_A=1$ is equivalent to $Q_B$ being a point.
This happens if and only if all vectors $b_i$ lie in a relevant line,
and
$\alpha_v=0$ for each relevant line $v$.
This last condition is equivalent to $P_B$ being
centrally symmetric.
\end{proof}

The following variant to the formula of Theorem \ref{algoth} works
well in practice for computing the $A$-discriminant $D_A$.
In the affine plane with coordinates $(w_1,w_2)$, consider
the following parametrically presented rational curve:
\begin{equation} \label{horn}
w_\ell \quad = \quad \prod_{i = 1}^n (b_{i1} + b_{i2} t)^{b_{i\ell}}
\,, \qquad \ell=1,2.
\end{equation}
This is the {\it Horn uniformization} in \cite[\S 9.3.C]{gkzbook}.
Let $\, \Delta(w_1,w_2) \,$ be the irreducible polynomial
defining this curve. This is a dehomogenization of the
$A$-discriminant, by Theorem \ref{algoth} or by
\cite[Theorem 9.3.3.~(a)]{gkzbook}.
More precisely,
\begin{equation}
\label{HornAlgo}
D_A(x_1,\ldots,x_n) \quad = \quad
(\hbox{a monomial}) \cdot \Delta(
\prod_{i=1}^n x_i^{b_{i1}},
\prod_{i=1}^n x_i^{b_{i2}}).
\end{equation}
The common factors in the
numerator and denominator of (\ref{horn})
are precisely the relevant lines which are not a
coordinate axis. In other words,
cancelling common factors in (\ref{horn})
is equivalent to replacing
$\,h_i(t) \, $ by $\,p_i(t) \,$ in (\ref{pes}).

We can get a description of the Newton polygon $N(\Delta)$ of
$\Delta(w_1,w_2)$ by ``dehomogenizing'' the result in
Theorem \ref{nq} as follows. Let
$b_i^\perp := (b_{i2}, - b_{i1})$ and note
that $- \det(b_i,v)$ equals the inner
product $ \langle b_i^\perp, v \rangle$.

\begin{corollary}
Let $B^\perp =\{b_1^\perp, \ldots, b_n^\perp\}$
and consider the polygon $Q_{B^\perp}$ translated
so that it lies in the first quadrant and
its boundary intersects both coordinate axes.
Then $N(\Delta) = Q_{B^\perp}.$
\end{corollary}

This result has been obtained independently
by Sadykov \cite{timur}, under the hypothesis that
there are no relevant lines outside
the coordinate axes.

\begin{example} \rm
We consider the toric $3$-fold of
degree $43$ in ${\bf P}^5$ which appears as Example $5.10$ in
\cite{ps}.
It is defined by the $6 \times 2$ integer matrix $B$
with row vectors $\,(2,3), \,
(-1,4),\, (-5,1), \, (3,-1), \, (2,-3),\,(3,-2) $.
The lattice ideal $I_B$ has seven minimal
generators. There are no relevant lines.
The polygon $P_B = Q_B$ is a hexagon.
Using Remark \ref{Pick} we find that $Q_B$ contains
$40$ lattice points. They correspond to the
$40$ terms in the $A$-discriminant $D_A$. The $6$ vertices
of $P_B$ correspond to the various leading terms in $D_A$.
Using (\ref{HornAlgo}) in any
computer algebra system we easily compute:
$$
\begin{array}{lcl}
D_A & = & -
\left (7\right )^{7}\left (17\right )^{17}\left (19\right )^{19} \,
{x_1}^{16} \,{x_4}^{11}\,{x_5}^{23}\,{x_6}^{22} \\
& & -\left (2\right )^{34}\left (3\right )^{15}\left (5\right
)^{15}\left(
13\right )^{13} \,
{x_1}^{20} \,{x_2}^{36}\, {x_3}^{11}\, {x_6}^{5}\,\\
& & \mbox{} + \left (2\right )^{10}\left (5\right )^{15}
\left (11\right )^{11}\left (17\right )^{17}
\,{x_1}^{23} \,{x_2}^{19} \, \,{x_5}^{13}{x_6}^{17} \\
& & \mbox{} + \left (2\right )^{64}\left (7\right )^{14}\left
(13\right )^{13}\,
{x_3}^{19} \,{x_4}^{28}\,{x_5}^{16} \,{x_6}^{9} \\
& & \mbox{} +\left (3\right )^{21}\left (7\right )^{7}\left
(11\right )^{11}
\left (13\right )^{13} \,
{x_2}^{16} \,{x_3}^{26}\,{x_4}^{25}\,{x_5}^{5} \\
& & \mbox{} -\left (2\right )^{10}\left (5\right )^{15} \left
(11\right )^{11}\left (17\right )^{17}
\,\,{x_1}^{9}\,{x_2}^{29} \, {x_3}^{21}\, {x_4}^{13} \\
& & \mbox{} + \, {\rm interior \, terms} .
\end{array}
$$
We invite the reader to draw $Q_B$ and verify
Theorem \ref{nq} for this example.
\qed
\end{example}

\section{Resultants having Newton triangles}

Mixed resultants form a subclass among all discriminants,
by the Cayley trick of elimination theory
\cite[\S 9.1.A]{gkzbook}. This subclass is important
for the theory of hypergeometric functions: conjecturally,
it consists of the denominators of rational hypergeometric functions
\cite[Conjecture 1.4]{rat}.
In this section we examine the Cayley construction and
mixed resultants in codimension $2$.

The Gale dual of a Cayley configuration $A$ is
a $(2r+3) \times 2$-matrix
$$ B \quad = \quad
\bigl( \, b_1,b_2,\ldots,b_r, \,c_1,c_2,\, -b_1,\ldots,-b_r,
\, -c_1-c_2 \, \bigr)^T,
$$
where the rows of the submatrix
$\, \tilde B \, := \,
( \, b_1,b_2,\ldots,b_r, \,c_1,c_2 \,)^T \,$
span $\Z^2$. We assume that all $b_i$ are non-zero
and ${\rm det}(c_1,c_2) \not= 0$.
By Corollary \ref{D1}, $D_A \not=1$.

Fix an $ r \times (r + 2)$-matrix
Gale dual to $\tilde B$ whose left
$r \times r$-minor is diagonal:
$$
\tilde A \,\, = \,\,
\left(
\begin{array}{ccccc}
\! \gamma_1 & & & \alpha_1 & \beta_1 \\
& \ddots & & \vdots & \vdots \\
& & \gamma_r & \alpha_r & \beta_r
\end{array} \right) \quad \hbox{where} \,\,
\gamma_i \in \Z_{> 0} \,\, \hbox{and} \,\,
(\alpha_i,\beta_i) \in \Z^2 \backslash \{(0,0)\}.
$$
This matrix lifts to a $(2r+1) \times (2r+3)$-matrix Gale dual to
$ B$ as follows:
\begin{equation}
A \quad = \quad
\left(
\begin{array}{cc}
\tilde A & 0 \\
I_{r+1} e_{r+1} & I_{r+1}
\end{array}
\right),
\end{equation}
where $I_{r+1}$ is the unit matrix of size $r+1$ and
$e_{r+1} = (0,0,\ldots,0,1)^{T}$.
The columns of $A$ index the coefficients in
a sparse system of $r+1$ equations:
\begin{eqnarray*}
f_0 \quad = &
z_1 \cdot t_1^{\alpha_1} \cdots t_r^{\alpha_r} \,+ \,
z_2 \cdot t_1^{\beta_1} \cdots t_r^{\beta_r} \,+ \,
z_3 \\
f_i \quad = & x_i \cdot t_i^{\gamma_i} \, + \, y_i \qquad
\hbox{for} \,\,\, i = 1,\ldots,r.
\end{eqnarray*}
This system consists of $r$ binomials and one Laurent trinomial,
as in (\ref{ThreeBinoOneTrino}).
The sparse resultant $\,{\rm Res}(f_0,f_1,\ldots,f_r)\,$ is the unique
(up to sign) irreducible polynomial in $
x_1,\ldots,x_r, y_1,\ldots,y_r,z_1,z_2,z_3$ which vanishes
when the system has a common root $(t_1,\ldots,t_r)$ in
the $r$-torus. {}From \cite[Prop. 9.1.7]{gkzbook} we get:

\begin{remark}
$\!\! $ The $A$-discriminant $D_A$ equals
the sparse resultant of $f_0,\ldots \!, \! f_r $.
\end{remark}

We now apply the product formula for resultants \cite{peders},
which amounts to evaluating $f_0$ at the common zeros of
$f_1,\ldots,f_r$.
The number of zeroes equals
$$ \Gamma \quad := \quad \gamma_1 \gamma_2 \cdots \gamma_r
\quad = \quad |\,{\rm det}(c_1,c_2) \,|. $$
Let $\eta_i$ denote a primitive $\gamma_i$-th root of unity.
The product formula implies:

\begin{proposition}
Up to a Laurent monomial factor, the $A$-discriminant is
$$
D_A \quad = \quad {\rm monomial} \cdot
\prod_{i_1=1}^{\gamma_1} \cdots
\prod_{i_r=1}^{\gamma_r}
f_0 \bigl( \eta_1^{i_1} z_1,
\eta_2^{i_2}z_2, \ldots,
\eta_r^{i_r} z_r
\bigr)
$$
where $z
_i = \left(- \frac{y_i}{x_i}\right)^{1/\gamma_i} \, $ for $\, i
=1,\ldots,r$.
\end{proposition}

Since $f_0$ is a trinomial, this formula gives an upper bound
of $\, \binom{\Gamma + 2}{2} \,$ for
the number of terms appearing in the expansion of $\,D_A =
{\rm Res}(f_0,\ldots,f_r)$.
This bound is quadratic in $\Gamma$.
In truth, this number grows linearly in $\Gamma$.

\begin{theorem}
The number of terms appearing in $D_A $ is at most
$\frac{5}{4} \cdot \Gamma \, + \, \frac{7}{4}$.
\end{theorem}

This bound is tight if the vectors $c_1 $ and $c_2$ span the lattice
$\Z^2$.
In this case, $\Gamma = {\rm det}(c_1,c_2) = 1$ and
the resultant $D_A$ has three terms. It is also
tight for the example in the Introduction, where
$\Gamma = 4$ and $D_A$ has six terms.

\begin{proof} According to Theorem \ref{nq},
the Newton polygon of the discriminant $D_A$ is essentially
the lattice triangle $\, Q_B \, = \, {\rm conv} \{ 0,c_1,c_2 \}$.
The number of terms in $D_A$ is at most the number of lattice points in
$Q_B $. Using Pick's formula as in Remark \ref{Pick}, we find
that the number $\, \# \,( Q_B \, \cap \, \Z^2) \,$ equals
$$
1 + \frac{1}{2} \cdot \bigl( \,|\det(c_1,c_2)|
\, + \,gcd(c_{11},c_{12})
\, + \,gcd(c_{21},c_{22})
\, + \,gcd(c_{11}+c_{21},c_{12}+c_{22}) \,\bigr).
$$
Using the inequality $\,a + b \leq ab + 1$, we
find that the sum of any two of the three last summands
is bounded above by $\, \Gamma + 1 = |\det(c_1,c_2)| + 1 $.
Therefore,
$$\, \# \,( Q_B \, \cap \, \Z^2) \quad \leq \quad
1 + \frac{1}{2} \cdot \bigl( \,\Gamma
\,+\, \frac{3}{2} \cdot (\Gamma +1) \,\bigr).
$$
This is the desired inequality.
\end{proof}

\bigskip
\bigskip

\noindent{\bf Acknowledgements:} We are grateful to Eduardo Cattani
for helpful discussions and to Laura Matusevich and
Thorsten Theobald for comments on
an earlier draft.
Alicia Dickenstein was supported by UBACYT TX94,
ANPCyT Grant 03-6568 and CONICET, Argentina.
Bernd Sturmfels was supported
by NSF Grant DMS-9970254 and the Miller Institute at UC Berkeley.

\smallskip

{\small
\begin{flushright}
Alicia Dickenstein \\
{\tt alidick@dm.uba.ar} \\
Dto.~de Matem\'atica, FCEyN, Universidad de Buenos Aires \\
(1428) Buenos Aires, \ Argentina
\end{flushright}

\begin{flushright}
Bernd Sturmfels \\
{\tt bernd@math.berkeley.edu} \\
Dept.~of Mathematics, University of California \\
Berkeley, CA 94720, USA
\end{flushright}}

\begin{thebibliography}{999}


\bibitem{rat} E.~Cattani, A.~Dickenstein, and B.~Sturmfels,
Rational hypergeometric functions, to appear in
{\it Compositio Mathematica}.

\bibitem{es} D.~Eisenbud and F.~Schreyer:
Chow forms and resultants via exterior algebra,
in preparation.


\bibitem{gkz1} I. M.~Gel'fand, A.~Zelevinsky, and M.~Kapranov:
{Hypergeometric functions and toric varieties},
{\it Funct.~Anal.~Appl.} {\bf 23} (1989) 94--106.

\bibitem{gkzbook} I. M.~Gel'fand, M.~Kapranov, and A.~Zelevinsky:
{\it Discriminants, Resultants and Multidimensional Determinants,}
Birkh\"auser, Boston, 1994.

\bibitem{mp} D.~Morrison and M.R.~Plesser: Summing the instantons,
quantum cohomology and mirror symmetry in toric varieties,
{\it Nuclear Physics B} {\bf 440} (1995) 279--354.

\bibitem{peders} P.~Pedersen and B.~Sturmfels:
Product formulas for resultants and Chow forms,
{\it Mathematische Zeitschrift} {\bf 214} (1993) 377--396.

\bibitem{ps} I.~Peeva and B.~Sturmfels:
Syzygies of codimension 2 lattice ideals,
{\it Mathematische Zeitschrift} {\bf 229} (1998) 163--194.

\bibitem{timur} T.~M.~Sadykov:
The Hadamard product of hypergeometric series,
Preprint, {\tt http://www.matematik.su.se/reports/2001/}.


\bibitem{StuSparse} B.~Sturmfels: Sparse elimination theory,
in {\sl ``Computational Algebraic Geometry and Commutative Algebra''}
[D.~Eisenbud and L.~Robbiano, eds.], Cambridge University Press,
1993, pp.~264--298.

\bibitem{Stubook} B.~Sturmfels: {\sl Gr\"obner Bases and
Convex Polytopes}, American Mathematical Society, 1995.


\end{thebibliography}
\end{document}